\documentclass[12pt,twoside,a4paper]{article}
\usepackage[cp866]{inputenc}
\usepackage[russian]{babel}
\usepackage{amsmath,amsfonts,amscd,amssymb,latexsym}

\begin{document}

\begin{center}

{\large\bf THE KANTOR--KOECHER--TITS CONSTRUCTION \hspace*{\fill} \\
FOR JORDAN COALGEBRAS} \hspace*{\fill} \\ [2mm]
\hspace*{6mm} {\large\bf V. N. Zhelyabin$^{*}$}
\footnote{$^{*}$Supported by ISF grant No. RB 6000.} \\
\hspace*{6mm}   E-mail vicnic@math.nsc.ru

\vspace{\baselineskip}
\end{center}

\begin{quote}\it{The relationship between Jordan and Lie coalgebras
is established. We prove that from any Jordan coalgebra $\langle
A, \Delta\rangle$, it is possible to construct a Lie coalgebra
$\langle L(A), \Delta_{L}\rangle$. Moreover, any dual algebra of
the coalgebra $\langle L(A), \Delta_{L}\rangle$ corresponds to a
Lie algebra that can be determined from the dual algebra for
$\langle A,\Delta\rangle$, following the Kantor--Koecher--Tits
process. The structure of subcoalgebras and coideals of the
coalgebra $\langle L(A), \Delta_{L}\rangle$ is
characterized.}\end{quote}

The notion of an (associative) coalgebra has long been traditionally
associated with the theory of Hopf algebras. The sharp increase of
interest in Hopf algebras was provoked by bringing into the focus
quantum groups, as well as the notion of a Lie bialgebra introduced
by Drinfeld in [1]. The definition of a coalgebra given in [2] is
related to a certain variety of algebras; in particular, a Jordan
coalgebra is defined as a coalgebra whose dual algebra is Jordan.

It is well known that to any Jordan algebra one can apply the
Kantor--Koecher--Tits (KKT) construction to produce a Lie algebra,
whose structure is closely connected with that of the initial Jordan
algebra. The natural question arises of whether an analog of the
KKT process is available for Jordan coalgebras. In this article we
provide an account that answers this question.

Let $k$ be a field. For linear spaces $U$ and $V$ over $k$,
$U\otimes V$ denotes their tensor product. By $V^*$ we denote the
space of all linear functionals given on $V$. For elements $f\in
V^*$ and $v\in V$, write $\langle f,v\rangle $ for a value of
the linear functional $f$ at $v$. If $X$ is a subset of $V$
$(V^*)$, then, as usual, $X^\bot$ is an orthogonal complement of
$X$ in the space $V^*$ $(V)$.

The pair $\langle A,{\Delta}\rangle $, where $A$ is a linear
space, and ${\Delta}\!:A\to A\otimes A$ is a linear mapping, is
called a {\it coalgebra.} The map ${\Delta}$, in this case, is
referred as {\it comultiplication.} For an element $a$ of the
coalgebra $\langle A,{\Delta}\rangle $, if ${\Delta}(a)=\sum
a_{1i}\otimes a_{2i}$, we write ${\Delta}(a)=\sum a_{(1)}\otimes
a_{(2)}$, following Sweedler [3].

Consider the dual space $A^*$. By the standard embedding
$\rho\!:A^*\otimes A^*\to (A\otimes A)^*$ we mean, as usual, a
linear map given by the following rule: $\langle \rho(f\otimes g),
\sum a_i\otimes b_i\rangle =\sum \langle f,a_i\rangle \langle
g,b_i\rangle$. Let ${\Delta}^*\!:(A\otimes A)^*\to A^*$ be the
dual of the map ${\Delta}$. On the space $A^*$, we can then
define multiplication by setting $fg={\Delta}^*\rho(f\otimes g)$
for any elements $f$ and $g$ in $A^*$. The space $A^*$ with
given multiplication is an ordinary algebra over a field $k$,
which we call the {\it dual algebra} of the coalgebra $A$. It is
easy to see that for any $f$, $g$ in $A^*$ and $a$ in $A$,
the equality $\langle fg,a\rangle=\langle \rho (f\otimes g),
{\Delta}(a)\rangle$ holds. This immediately implies the following:
$\langle fg,a\rangle=\sum \langle f,a_{(1)}\rangle \langle
g,a_{(2)}\rangle$, where ${\Delta}(a)=\sum a_{(1)}\otimes a_{(2)}$.

The linear subspace $B$ of $\langle A,{\Delta}\rangle$ is called a
{\it subcoalgebra} ({\it coideal}) if ${\Delta} (B)\subseteq
B\otimes B$ $({\Delta}(B)\subseteq B\otimes A+A\otimes B)$. The
coalgebra $\langle A,{\Delta}\rangle$ is said to be {\it locally
finite-dimensional} if each finitely generated subcoalgebra is
finite-dimensional. If $A^*$ is an algebra with unity 1, we assume
that $1\in B^\bot$ for any coideal $B$ of $\langle
A,{\Delta}\rangle $. We know from [2] that $B$ is a subcoalgebra
(coideal) in $\langle A,{\Delta}\rangle $ if and only if $B^\bot$
is an ideal (subalgebra) of $A^*$. The dual algebra $A^*$ of
$\langle A,{\Delta}\rangle $ induces a bimodule action $(\cdot )$
on the space $A$, defined as follows:
$$ f\cdot a=\sum a_{(1)} \langle f,a_{(2)}\rangle\mbox{\, and \,}
a\cdot f=\sum \langle f,a_{(1)}\rangle a_{(2)}, $$

\noindent where $f\in A^*$ and ${\Delta}(a)=\sum a_{(1)}\otimes
a_{(2)}$. It is easy to see that for any $f$, $g\in A^*$ and
$a\in A$, $\;\langle fg,a\rangle=\langle f, g\cdot a\rangle =
\langle g, a\cdot f\rangle$. It is not hard to show that $B$ is
a subcoalgebra in $\langle A,{\Delta}\rangle$ if and only if $B$
is an $A^*$-subbimodule of the bimodule of $A$.

Now let $M$ be some variety of algebras over a field $k$. Then
$\langle A,{\Delta}\rangle$ will be called an $M$-{\it coalgebra}
if its dual algebra $A^*$ belongs to $M$. Specifically, if the
dual algebra $A^*$ is Jordan, then $\langle A,{\Delta}\rangle $ is
a {\it Jordan coalgebra}. In [2], it was shown that every Jordan
coalgebra is locally finite-dimensional and is defined as a Jordan
bimodule over its dual algebra.

\medskip
\begin{center}
{\bf 1. WEAK INNER DERIVATIONS}
\end{center}
\smallskip

Let $J$ be a Jordan algebra over a field $k$ (of characteristic
not 2). The map ${\rm d}\!:J\to J$ is called a {\it derivation} if
$(ab){\rm d}=a{\rm d}b +b{\rm d}a$ for any elements $a$ and $b$
in $J$. The linear space ${\rm Der}(J)$ of all derivations of
$J$ is a Lie algebra w.r.t. commutation. If $a$ is an element of
$J$, $a'$ stands for the operator of right multiplication by $a$.
Let $X$ be a subset of $J$. We then put $X'=\{a'|a\in X\}$.
For elements $a$ and $b$ in $J$, write $[a,b]=a'b'- b'a'$.
For subsets $X$ and $Y$ in $J$, write $[X,Y]$ to denote the
linear space generated by elements $[x,y]$, where $x\in X$ and
$y\in Y$. It is easy to see that for any $a$, $b$, $c$,
and $d$ in $J$, the following relations hold:
\begin{equation}
[ab,c]+[ac,b]+[bc,a]=0
\end{equation}

\noindent and
\begin{equation}
(cd)[a,b]=c[a,b]d'+d[a,b]c'.
\end{equation}

\noindent By (2), the operator $[a,b]$ is a derivation of $J$.
Derivations of $J$ of the form $\sum [a_i,b_i]$ are said to be
{\it inner}. Denote by ${\rm IntDer} (J)$ the linear space
generated by all inner derivations. The space ${\rm IntDer} (J)$ is
then an ideal in the Lie algebra ${\rm Der}(J)$. Suppose that the
characteristic of the field $k$ is not 2; $\langle A,
{\Delta}\rangle$ is a Jordan coalgebra and $A^*$ is its dual
algebra.

\smallskip
{\bf Definition.} A derivation ${\rm d}$ of $A^*$ is called {\it
weakly inner} if, for any finite-dimensional subcoalgebra $B$ in
$\langle A,{\Delta}\rangle $, there exists a derivation $\rm i$
from ${\rm IntDer}(A^*)$ such that $\langle f{\rm d}, b\rangle
=\langle f{\rm i}, b\rangle $ for any $f\in A^*$ and $b\in B$.

Denote by ${\rm WIntDer} (A^*)$ the linear space of all weakly
inner derivations of $A^*$. Obviously, ${\rm
IntDer}(A^*)\subseteq {\rm WIntDer}(A^*)$. It is easy to see that
a derivation ${\rm d}$ of $A^*$ belongs to ${\rm WIntDer}(A^*)$
if and only if, for any finite-dimensional subcoalgebra $B$ in
$\langle A,{\Delta}\rangle$, there is a derivation $\rm i$ of
${\rm IntDer} (A^*)$ such that $A^*({\rm {d-i}})\subseteq B^\bot$.

Suppose $\rm i$ is an inner derivation of $A^*$. Then $\rm i$
naturally induces some linear map of the space $A$ into itself,
which we denote also by $\rm i$. If ${\rm i}=\sum [f_i,g_i]$ and
$a\in A$, put $a{\rm i}=\sum (f_i,a,g_i)$. It is not hard to see
that for any elements $x$, $y$, $z\in A^*$ and $a\in A$, the
equality $$\langle (x,y,z),a\rangle =-\langle y,(x,a,z)\rangle $$
holds. Therefore, for any elements $f$ in $A^*$ and $a$ in $A$, we
obtain $\langle f{\rm i},a\rangle =-\langle f,a{\rm i}\rangle$.

\smallskip
{\bf LEMMA 1.} The linear space ${\rm WIntDer}(A^*)$ is a
subalgebra of the Lie algebra ${\rm Der}(A^*)$.

{\bf Proof.} Suppose that elements ${\rm d}_1$ and ${\rm d}_2$
belong to ${\rm WIntDer}(A^*)$. In virtue of the definition of
${\rm WIntDer}(A^*)$, then, for any finite-dimensional
subcoalgebra $B$ in $A$, there are ${\rm i}_1$ and ${\rm i}_2$
from ${\rm IntDer}(A^*)$ such that for any $f\in A^*$ and $b\in
B$ we have $\langle f{\rm d}_1,b\rangle =\langle f{\rm
i}_1,b\rangle$ and $\langle f{\rm d}_2,b\rangle =\langle f{\rm
i}_2,b\rangle $. Now consider $\langle f({\rm d}_1{\rm d}_2- {\rm
d}_2{\rm d}_1),b\rangle$. We can show that $\langle f({\rm
d}_1{\rm d}_2-{\rm d}_2{\rm d}_1),b\rangle = - \langle f{\rm
d}_1,b{\rm i}_2\rangle + \langle f{\rm d}_2,b{\rm i}_1\rangle$.
Since $B$ is a subcoalgebra, $b{\rm i}_1$, $b{\rm i}_2\in B$.
Therefore,
$$ \langle f({\rm d}_1{\rm d}_1-{\rm d}_2{\rm d}_1),b\rangle=- \langle f{\rm i}_1,b{\rm i}_2\rangle +\langle f{\rm i}_1,b{\rm i}_2\rangle = \langle f({\rm i}_1{\rm i}_2-{\rm i}_2{\rm i}_1), b\rangle. $$

\noindent But ${\rm i}_1{\rm i}_2-{\rm i}_2{\rm i}_1$ is an inner
derivation, whence ${\rm d}_1{\rm d}_2- {\rm d}_2{\rm d}_1 \in {\rm
WIntDer}(A^*)$. Thus, ${\rm WIntDer}(A^*)$ is a subalgebra of the
Lie algebra ${\rm Der}(A^*)$, as was to be proved.

Below we assume that the algebra $A^*$ contains unity, and that
this is equivalent to the possession of counity by the coalgebra
$\langle A,{\Delta}\rangle$. In the space of linear endomorphisms
of $A^*$, consider a subspace $R(A^*)= (A^*)' \oplus {\rm
WIntDer}(A^*)$. The space $R(A^*)$ is a Lie algebra with respect
to commutation $[\,,\,]$. Denote by ${\varepsilon}$ the
mapping from $R(A^*)$ into $R(A^*)$, given by the formula
${\varepsilon}(f'+{\rm d})= -f'+{\rm d}$. Consider the linear
space $L(A^*)=A^*\oplus R(A^*) \oplus \bar{A}^*$, where
$\bar{A}^*$ is an isomorphic copy of the space $A^*$. Define
multiplication $[\,,\,]$ on $L(A^*)$ by setting
$$ [a_1+{\rm d}_1+\bar{b}_1, a_2+{\rm d}_2+\bar{b}_2]= a_1{\rm d}_2-a_2{\rm d}_1+a_1\nabla b_2- a_2\nabla b_1 + [{\rm d}_1,{\rm d}_2] + \overline{b_1{\varepsilon}({\rm d}_2)} - \overline{b_2{\varepsilon}({\rm d}_1)}, $$

\noindent where $a_1$, $a_2$, $b_1$, $b_2\in A$, ${\rm d}_1$,
${\rm d}_2\in R(A^*)$, and $a\nabla b=(ab)'-[a,b]$. We know
that the algebra thus obtained is a Lie algebra. The space
$I(A^*)=A^*\oplus (A^*)'\oplus {\rm IntDer}(A^*) \oplus
\bar{A}^*$ is an ideal in $L(A^*)$. Extend the map
${\varepsilon}$ to the whole algebra $L(A^*)$ by setting
${\varepsilon}(a+c'+{\rm d}+\bar{b})=b-c'+{\rm d}+\bar{a}$.

Consider the algebra $C=A^*+A$, the null extension of $A^*$. Let
$X$ be a subset of $A$. We can then put ${\rm Ann} (X)=\{f\in
A^*|f\cdot X =0\}$. Recall that for subsets $X$ and $Y$ in $C$,
$[X,Y]$ is the linear space generated by operators $[x,y]$,
where $x\in X$ and $y\in Y$.

\smallskip
{\bf Proposition 1.} Let $B$ be a subcoalgebra of $(A,{\Delta})$.
Then ${\rm Ann}(B)=B^\bot$ and $[B^\bot,B]=0$. If $B$ is a
finite-dimensional space, the space $[A^*, B]$ is also
finite-dimensional.

{\bf Proof.} Since $B^\bot$ is an ideal in $A^*$ (see [2]), we
have $\langle f,B^\bot\cdot B\rangle=\langle fB^\bot, B\rangle =0$.
Consequently, $B^\bot\cdot B =0$, and so $B^\bot\subseteq {\rm
Ann}(B)$. Since $\langle {\rm Ann}(B),B\rangle =\langle 1,{\rm
Ann}(B)\cdot B\rangle =0$, it follows that ${\rm Ann}(B)\subseteq
B^\bot$. Thus, $B^\bot = {\rm Ann}(B)$. Hence, $C[B^{\bot},
B]\subseteq B^\bot \cdot B = 0$, that is, $[B^\bot, B]=0$. Let
$B$ be a finite-dimensional space. Then $B^\bot$ has a finite
codimension. Therefore, $A^*=V+B^\bot$, where $V$ is a
finite-dimensional subspace, whence $[A^*,B]=[V,B]$. Thus, the
space $[A^*,B]$ is finite-dimensional, as desired.

\smallskip
{\bf Proposition 2.} Every derivation ${\rm d}$ of ${\rm
WIntDer}(A^*)$ can be extended to a derivation $D$ in $C$, in
which case for any subcoalgebra $B$ of $\langle A,{\Delta}\rangle$,
the inclusion $BD\subseteq B$ holds, and we have $\langle f{\rm
d},a\rangle =-\langle f,aD\rangle$ for any $f\in A^*$ and $a\in A$.
If the subcoalgebra $B$ is finite-dimensional, there exists an
element ${\rm i}\in {\rm IntDer}(A^*)$ such that $bD=b{\rm i}$ and
$\langle fD,b\rangle=\langle f{\rm i},b\rangle$ for any $b\in B$
and $f\in A^*$.

{\bf Proof.} Let $d\in {\rm WIntDer}(A^*)$, $a\in A$, and $B$
be a finite-dimensional subcoalgebra, containing $a$. There then
exists ${\rm i}\in {\rm IntDer}(A^*)$ such that $\langle f{\rm d},
b\rangle =\langle f{\rm i},b\rangle $ for any $f\in A^*$ and $b\in
B$. Since $a\in B$, we put $aD=a{\rm i}$ and $fD=f{\rm d}$ for
all $f\in A^*$. The map $D$ does not depend on the choice of $B$.
Indeed, let $B_1$ be another finite-dimensional subcoalgebra
containing $a$. Then there exists ${\rm i}_1\in {\rm IntDer}(A^*)$
such that $\langle f{\rm d},b\rangle =\langle f{\rm i}_1,b\rangle$
for any $f\in A^*$ and $b\in B_1$. Therefore, for all $f$ in
$A^*$ we have $\langle f{\rm d},a\rangle =\langle f {\rm i}_1,
a\rangle $. This implies that $\langle f,a({\rm i-i_1})\rangle = -
\langle f({\rm i-i_1}),a\rangle =0$ for any $f\in A^*$, that is,
$a{\rm i}=a{\rm i}_1$.

Obviously, $\langle f{\rm d},b\rangle =-\langle f,bD\rangle $ for any $f\in A^*$
and $b\in A$, and $BD\subseteq B$ for any subcoalgebra $B$ in
$\langle A,{\Delta}\rangle$. Suppose that a subalgebra $B$ in $A$
is finite-dimensional. There then exists an element ${\rm i}\in
{\rm IntDer}(A^*)$ such that $\langle f{\rm d},b\rangle= \langle
f{\rm i},b\rangle$ for any $f\in A^*$ and $b\in B$. Therefore,
$\langle f,bD\rangle =-\langle f{\rm d},b \rangle= \langle f,b{\rm
i}\rangle$, that is, $bD=b{\rm i}$ for any $b$ in $B$. We show
that $D$ is a derivation of the algebra $C$. Let $f$, $g\in
A^*$ and $a\in A$. Since $g{\rm d}=gD$ and $\langle f,(a\cdot
g)D\rangle =-\langle f{\rm d},(a\cdot g)\rangle $, we have
$$ \langle f,(a\cdot g)D\rangle=-\langle f{\rm d}g,a\rangle = -
\langle (fg){\rm d},a\rangle + \langle g{\rm d}f, a\rangle =
\langle f, aD\cdot g\rangle + \langle f,gD\cdot a\rangle.$$

\noindent Hence, $(a\cdot g)D=aD\cdot g+gD\cdot a$, and so $D$ is
a derivation of $C$, as desired.

\smallskip
{\bf LEMMA 2.} There exists an isomorphism $\phi$ between the space
${\rm WIntDer}(A^*)$ and the space of all linear functionals,
defined on the space $[A^*,A]$, such that $\langle \phi ({\rm d}),
[f,a]\rangle =\langle f{\rm d},a\rangle $ for any $f\in A^*$,
$a\in A$, and ${\rm d}\in {\rm WIntDer}(A^*)$.

We divide the proof of the lemma into a series of propositions.

\smallskip
{\bf Proposition 3.} There exists an isomorphic embedding $\phi$ of
the space ${\rm WIntDer}(A^*)$ into the space of all linear
functionals, given on the space $[A^*,A]$, such that $\langle
\phi({\rm d}), [f,a]\rangle =\langle f{\rm d},a\rangle$ for any
$f\in A^*$, $a\in A$, and $d\in {\rm WIntDer}(A^*)$. If the
coalgebra $\langle A,{\Delta}\rangle $ is finite-dimensional,
then $\phi$ is an isomorphism.

{\bf Proof.} Let ${\rm d}$ be an element in ${\rm WIntDer}(A^*)$
and $v$ an arbitrary element from $[A^*,A]$. If $v=\sum\limits_i
[f_i, a_i]$, the functional $\hat{\rm d}$ on $[A^*,A]$ is then
defined by setting $\langle \hat{\rm d},v\rangle =\sum\limits_i
\langle f_i{\rm d},a_i\rangle$. We show that $\hat{\rm d}$ is
well defined. Indeed, let $u$ be an element in $[A^*,A]$ such
that $u=v$. Assume $u=\sum\limits_j [g_j,b_j]$. Suppose that
$B$ is a subcoalgebra of $(A,{\Delta})$ generated by all
elements $a_i$ and $b_j$. Then $B$ is finite-dimensional. Since
${\rm d}\in {\rm WIntDer}(A^*)$, there exists an element ${\rm i}$
in ${\rm IntDer}(A^*)$ such that $\langle f{\rm d},b\rangle=
\langle f{\rm i},b\rangle$ for any $f\in A^*$ and $b\in B$.
Since ${\rm i}=\sum\limits_k [h_k,e_k]$, where $h_k$, $e_k\in
A^*$, we have
$$ \langle \hat{\rm d},v\rangle = \sum\limits_i \langle f_i{\rm d},a_i\rangle =\sum\limits_i\sum\limits_k \langle f_i[h_k,e_k], a_i\rangle = \sum\limits_i\sum\limits_k \langle e_k, h_k [f_i, a_i] \rangle = $$
$$ \sum\limits_k \langle e_k, h_k v\rangle =\sum\limits_k \langle e_k,h_ku\rangle = \sum\limits_j\sum\limits_k \langle g_j[h_k,e_k], b_j\rangle = \sum\limits_j \langle g_j{\rm d},b_j\rangle = \langle \hat{\rm d},u\rangle. $$

\noindent This means that the functional $\hat{\rm d}$ is well
defined. It is now evident that the map $\phi\!:{\rm
WIntDer}(A^*)\to [A^*,A]^*$, given by $\phi ({\rm d})=\hat{\rm d}$,
is the desired isomorphic embedding.

Let $\langle A,{\Delta}\rangle$ be a finite-dimensional coalgebra.
Clearly, ${\rm WIntDer}(A^*)={\rm IntDer}(A^*)$.
Consider a nonzero element $u$ from $[A^*,A]$. Since $(A^*)u\neq
0$, it follows that $fu\neq 0$ for some element $f$ in $A^*$.
Consequently, $\langle g, fu\rangle\neq 0$ for some functional $g$
from $A^*$, and it is easy to see that $\langle \phi
([f,g]),u\rangle =\langle g,fu\rangle$. Therefore, $u$ determines
the nonzero functional on the space $[A^*,A^*]$. Hence, $\phi$
is an isomorphism, as desired.

Now we let $B$ be a finite-dimensional subcoalgebra of $\langle
A,{\Delta}\rangle$. Then $A^*=V\oplus B^\bot$, where $V$ is a
finite-dimensional subspace in $A^*$.

\smallskip
{\bf Proposition 4.} The space $V$ can be given the structure of an
algebra with multiplication $*$, so that the resulting algebra $V$
is the dual of the coalgebra $\langle B,{\Delta}\rangle$, and for
any elements $f$, $g\in V$ and $b\in B$, the equality $\langle
f*g,b\rangle =\langle fg,b\rangle $ holds.

{\bf Proof.} Let $\rho$ be the standard embedding of the space
$A^*\otimes A^*$ into $(A\otimes A)^*$. Obviously, $A=B\oplus
W$ and $V=B^*$. For the elements $v$ in $A^*\otimes A^*$ and
$f$ in $A^*$, denote by $\rho (v)_0$ and $f_V$ their
projections onto the spaces $(B\otimes B)^*$ and $V$,
respectively. Define the map $\rho_{B}\!:B^*\otimes B^*\to (B\otimes
B)^*$ by setting $\rho_B(f\otimes g)=\rho (f\otimes g)_0$.
Clearly, $\rho_B$ is a linear map. Since ${\Delta}$ is a linear
map from $B$ into $B\otimes B$, it possesses the dual map
${\Delta}^*$ from $(B\otimes B)^*$ into $B^*$, which is also
linear. And ${\Delta}^*\rho_B$ is then a linear map from $V\otimes
V$ into $V$.

First, let $f$ and $g$ be in $V$. We can then put
$f*g={\Delta}^*\rho_B (f\otimes g)$. It is not hard to see that
the space $V$ with $*$ is the dual algebra of the coalgebra
$\langle B,{\Delta}\rangle$. Moreover, for any element $b$ in $B$,
we have
$$ \langle f*g,b\rangle =\langle {\Delta}^*\rho_B(f\otimes g), b\rangle = \langle \rho_B (f\otimes g), {\Delta}(b)\rangle = \langle \rho(f\otimes g)_0, {\Delta}(b)\rangle. $$

\noindent Since $\rho(V\otimes B^\bot + B^\bot\otimes V +
B^\bot\otimes B^\bot)\subseteq (B\otimes B)^\bot$, it follows that
$\langle \rho (f\otimes g)_0, {\Delta}(b)\rangle = \langle \rho
(f\otimes g), {\Delta}(b)\rangle$. Consequently, $\langle
f*g,b\rangle =\langle fg,b\rangle$.

Next, let $D=V+B$ be the null extension of $V$. Denote
multiplication in the algebra $D$ by $*$. For elements $a$ and
$b$ in $D$, put $\{ a,b\}=a'b'- b'a'$, where $a'$ and $b'$
are operators of right multiplication in $D$. If $X$ and $Y$
are subsets of $D$, $\{ X,Y\}$ is the linear space generated by
operators $\{ x,y\}$, where $x\in X$ and $y\in Y$. It is easy
to see, then, that $f*a=f\cdot a$ for any elements $f$ in $V$ and
$a$ in $B$.

Further, let $f\in V$ and $b\in B$. Then, for any element $g$ in
$V$, the equality $g\{ f,b\} =g[f,b]$ holds. Indeed, $g\{ f,b\}
=(g*f)*b- g*(f*b)=\sum b_{(1)}\langle g*f,b_{(2)}\rangle - g\cdot
(f\cdot b)=\sum b_{(1)} \langle gf,b_{(2)}\rangle -g\cdot (f\cdot
b)= (gf)\cdot b-g\cdot (f\cdot b)=g[f,b]$. Hence,
$\sum\limits_{i=1}\{ f_{i}, b_{i}\}=0$ if and only if
$\sum\limits_{i=1}[f_{i},b_{i}]=0$, where $f_{i}\in V$,
$b_{i}\in B$, and $i=1, \ldots, n$. In this case, however, it
is not hard to see that the spaces $\{ V,B\}$ and $[V,B]$ are
isomorphic. By Proposition 3, we can therefore assume that $\{
V,V\}$ is the space of all linear functionals on $[V,B]$. In
addition, for any elements $v$ in $\{ V,V\}$, $f$ in $V$, and
$b$ in $B$, we have $\langle v,[f,b]\rangle =\langle fv,b\rangle$.

Suppose $v\in \{ V,V\}$ and $v\neq 0$. Then the functional $v$
can be extended to the whole space $[A^{*},A]$. Consider a
nonzero element $u$ in $[V,V]$. By Proposition 3, we can assume
that $u\in [A^{*},A]^{*}$. Consequently, for any $f\in V$ and
$b\in B$ we have $\langle v-u,[f,b]\rangle =\langle
f(v-u),b\rangle$. The elements $v$ and $u$ can be chosen in such
a way that $f(v-u)=0$. Therefore, $[V,V]\subseteq
[V,B]^{*}+[V,B]^{\perp}$. We thus arrive at the following:

\smallskip
{\bf Proposition 5.} Let $B$ be a finite-dimensional subcoalgebra
of the coalgebra $\langle A,\bigtriangleup\rangle$ and let $V$ be
the subspace in $A^{*}$ such that $A^{*}=V\oplus B^{\perp}$.
Then the space $\{ V,B\}$ is isomorphic to $[V,B]$, and the
following inclusion holds: $[V,V]\subseteq [V,B]^{*}+[V,B]^{\perp}$,
where $[V,B]^{\perp}$ is an orthogonal complement of $[V,B]$ in
the space $[A^{*},A]^{*}$.

We are now in a position to proceed to the

{\bf Proof} of Lemma 5. Let $\hat{\rm d}$ be a linear functional on
$[A^{*},A]$. Then we define the map ${\rm d}:\! A^{*}\longrightarrow
A^{*}$ by setting $\langle f{\rm d},a\rangle = \langle \hat{\rm
d},[f,a]\rangle$ for any $f$ in $A^{*}$ and $a$ in $A$. We show that
${\rm d}\in {\rm WIntDer}(A^*)$. By the definition of ${\rm d}$, in
view of (1) we have $$ \langle (fg){\rm d},a\rangle=\langle \hat{\rm
d},[fg,a]\rangle = \langle \hat{\rm d}, [f,ga]\rangle + \langle
\hat{\rm d}, [g,fa]\rangle = \langle f{\rm d},g\cdot a\rangle
+ \langle g{\rm d}, f\cdot a\rangle = \langle f{\rm d}g+g{\rm d}f,
a\rangle. $$

\noindent Consequently, ${\rm d}$ is a derivation of $A^*$. Let
$B$ be a finite-dimensional subcoalgebra of $\langle
A,{\Delta}\rangle $ such that $\langle (A^*){\rm d},B\rangle \neq
0$ and let $U= [A^*,B]$. Then $A^*=V\oplus B^\bot$, where $V$
is a finite-dimensional subspace, and $U=[V,B]$ by Proposition 1.
Since $\langle (A^*){\rm d},B\rangle\neq 0$, we have $\langle
\hat{\rm d},U\rangle\neq 0$.

First consider the space $\{V,B\}$, introduced above. By
Proposition 5, $\hat{\rm d}$ can be assumed to be a nonzero
functional, given on the space $\{V,B\}$. In addition, for any
elements $f$ in $V$ and $b$ in $B$, the equalities $\langle
\hat{\rm d}, \{f,b\}\rangle=\langle \hat{\rm d},[f,b]\rangle=
\langle f{\rm d},b\rangle$ hold. By Proposition 3, there then
exists an inner derivation $\sum\{f_i,g_i\}$ of $V$ such that
$\langle\hat{\rm d}, \{f,b\}\rangle=\langle f\sum\{f_i,g_i\},b
\rangle$.  By Proposition 4, $\langle f\sum\{f_i,g_i\},b\rangle =
\langle f\sum[f_i,g_i], b\rangle$. Therefore, $\langle f{\rm
d},b\rangle =\langle \hat{\rm d},[f,b]\rangle =\langle f\sum
[f_i,g_i],b\rangle$, that is, ${\rm d}\in {\rm WIntDer}(A^*)$. It is
now evident that the map $\phi\!: {\rm WIntDer}(A^*)\to [A^*,A^*]$,
given by $\phi({\rm d})= \hat{\rm d}$, is the desired isomorphism.

We continue to consider the space $L(A)=A\oplus R(A)\oplus\bar{A}$,
where $\bar{A}$ is a an isomorphic copy of the space $A$ and
$R(A)=A'\oplus [A^*,A]$. By Lemma 2, we can assume that ${\rm
WIntDer}(A^*)=[A^{*},A]^*$. Then the following is valid:

\smallskip
{\bf COROLLARY.} $L(A^*)$ is the space of all linear functionals
given on the space $L(A)$.

By Proposition 2, every derivation ${\rm d}$ from ${\rm
WIntDer}(A^*)$ can be extended to a derivation in the algebra $C$.
Denote the extended derivation also by ${\rm d}$. The Lie
algebra ${\rm Hom}_k(C,C)^{(-)}$ then satisfies the
equalities $[a',{\rm d}] =(a{\rm d})'$ and $[[f,a],{\rm
d}]=[[f{\rm d}, a]+[f,a{\rm d}]$, where $f\in A^*$ and $a\in A$.
Let $\pi$ be the map from $L(A)$ into $L(A)$, defined by
${\pi}(a+b'+v+\bar{c})= c-b'+v+\bar{a}$. Also let $a$, $c\in A$,
$\bar{b}\in \bar{A}$, $w\in A'\oplus [A^*,A]$, $f\in A^*$,
$\bar{g}\in (\bar{A}^*)$, $u\in (A^*)'\oplus {\rm WIntDer}(A^*)$,
and $a\nabla g=(ag)'+[g,a]$. On the space $L(A)$, define the
structure of an $L(A^*)$-module by setting
$$ [a+w+\bar{b}, f+u+\bar{g}]=au-fw+a\nabla g-f\nabla b +[w,u]+ \overline{b{\varepsilon}(u)} - \overline{g\pi (w)}. $$

\noindent It is well known that the action of the Lie algebra
$L(A^*)$ on $L(A)$, thus defined, transforms the space $L(A)$
into a Lie $L(A^*)$-bimodule.

\smallskip
{\bf LEMMA 3.} Let $l_1^*$, $l_2^*$ be elements in $L(A^*)$ and
let $l\in L(A)$. Then $\langle [l_1^*,l_2^*]$, $l\rangle
=\langle l_2^*$, $[l,{\varepsilon}(l_1^*)]\rangle $.

{\bf Proof.} Suppose $l_1^*=f_1+g'_1+{\rm d}_1+\bar{h}_1$,
$l_2^*=f_2+g'_2+{\rm d}_2+\bar{h}_2$, and $l=a+b'+v+\bar{c}$,
where $f_i$, $g_i\in A^*$, $d_i\in {\rm WIntDer}(A^*)$,
$\bar{h}_i\in (A^*)$ and $a$, $b\in A$, $v\in [A^*,A]$,
$\bar{c}\in \bar{A}$. Then
$$ [l_1^*,l_2^*]=f_1(g'_2+{\rm d}_2)-f_2(g'_1{\rm d}_1)+ f\nabla h_2 -f_2h_1 + $$
$$ [g'_1+{\rm d}_1, g'_2+{\rm d}_2]-
\overline{h_1g_2}+\overline{h_2g_1} +\overline{h_1{\rm d}_2} -
\overline{h_2{\rm d}_1}. $$

\noindent Hence, $\langle [l_1^*,l_2^*],l\rangle =\langle
f_1(g'_2+{\rm d}_2) - f_2(g'_1+{\rm d}_1),a\rangle + \langle
-\overline{h_1g_2}+ \overline{h_2g_1}+\overline{h_1{\rm
d}_2}-\overline{h_2{\rm d}_1}, \bar{c}\rangle + \langle f_1\nabla
h_2-f_2\nabla h_1 + [g'_1+{\rm d}_1, g'_1{\rm d}_2], b'+v\rangle$.
By Lemma 2 and Proposition 2, we have $\langle f_1(g'_2+{\rm
d}_2)-f_2(g'_1+{\rm d}_1)$, $a\rangle = \langle g_2$, $f_1\cdot
a\rangle +\langle {\rm d}_2$, $[f_1,a]\rangle +\langle f_2$,
$a{\rm d}_1\rangle$. Similarly, $\langle
-\overline{h_1g_2}+\overline{h_2g_1}+\overline{h_1{\rm d}_2}
-\overline{h_2{\rm d}_1}$, $\bar{c}\rangle =- \langle g_2$,
$h_1\cdot c\rangle +\langle h_2 $, $g_1\cdot c\rangle +\langle
{\rm d}_2$, $[h_1,c]\rangle +\langle h_2, c{\rm d}_1\rangle$.

Let $v=\sum[e_i,a_i]$. By Lemma 2, we obtain
$$ \langle -[f_1,h_2]+[f_2,h_1]+[g_1,g_2]+ [{\rm d}_1,{\rm d}_2], v\rangle= $$
$$ \sum\langle -[f_1,h_2]+[f_2,h_1]+[g_1,g_2]+[{\rm d}_1,{\rm d}_2],
[e_i,a_i]\rangle = $$
$$ \sum \langle -(f_1,e_i,h_2)+(f_2,e_i,h_1)+ (g_1,e_i,g_2),a\rangle
+ \sum\langle e_i[{\rm d}_1,{\rm d}_2],a_i\rangle. $$

\noindent Again, $$\sum\langle -(f_1,e_i,h_2)+ (f_2,e_i,h_1) +
(g_1,e_i,g_2), a_i\rangle =$$
$$\sum -\langle h_2,(e_i,f_1,a_i)\rangle - \langle f_2, (e_i,h_1,a_i)\rangle + \langle g_2,(e_i,g_1,a_i)\rangle = $$
$$ - \langle h_2,f_1 v\rangle - \langle f_2,h_1v\rangle + \langle
g_2,g_1v\rangle. $$

\noindent Lemma 2 and Proposition 2 yield
$$ \sum\langle e_i[{\rm d}_1,{\rm d}_2],a_i\rangle = \sum \langle
{\rm d}_2, [e_i{\rm d}_1,a_i]+[e_i,a_i{\rm d}_1]\rangle =
\sum\limits \langle {\rm d}_2, [[e_i,a_i], d_1]\rangle = \langle
{\rm d}_2,[v,d_1]\rangle, $$

\noindent from which we obtain
$$ \langle f_1\nabla h_2- f_2\nabla h_1 +[g'_1+{\rm d}_1, g'_2+ {\rm d}_2], b'+v\rangle =\langle (f_1h_2)'-(f_2h_1)',b'\rangle+ $$
$$ \langle (g_1{\rm d}_2)'-(g_2{\rm d}_1)', b'\rangle + \langle - [f_1,h_2]+[f_2,h_1]+[g_1,g_2]+[{\rm d}_1,{\rm d}_2], v\rangle = $$
$$ \langle f_1h_2-f_2h_1+g_1{\rm d}_2-g_2{\rm d}_1,b\rangle - \langle h_2,f_1v\rangle - \langle f_2,h_i v\rangle + \langle g_2,g_1v\rangle + \langle {\rm d}_2, [v,{\rm d}_1]\rangle. $$

\noindent Since $\langle f_1h_2-f_2h_1+g_1{\rm d}_2-g_2{\rm
d}_1,b\rangle =\langle h_2,f_1\cdot b\rangle - \langle f_2,h_1\cdot
b\rangle + \langle d_2,[g_1,b]\rangle + \langle g_2,b{\rm
d}_1\rangle$, we have $\langle f_1\nabla h_2- f_2\nabla h_1 +
[g'_1+{\rm d}_1, g'_2+{\rm d}_2], b'+v\rangle=\langle h_2,f_1\cdot
b\rangle - \langle f_2, h_1\cdot b\rangle +\langle {\rm d}_2,
[g_1,b]\rangle + \langle g_2, b{\rm d}_1\rangle -\langle
h_2,f_1v\rangle - \langle f_2, h_1v\rangle+\langle
g_2,g_1v\rangle+\langle {\rm d}_2,[v, {\rm d}_1]\rangle$. Thus
$\langle [l_1^*,l_2^*],l\rangle = \langle f_2, -g_1\cdot a+ a{\rm
d}_1-h_1\cdot b -h_1\cdot v\rangle+ \langle g_2,f_1\cdot a-h_1\cdot
c+g_1v +b{\rm d}_1\rangle +\langle {\rm d}_2,
[f_1,a]+[h_1,c]+[g_1,b]+[v,{\rm d}_1]\rangle + \langle h_2, g_1\cdot
c+c{\rm d}_1+f_1\cdot b- f_1v\rangle$. It is easy to verify that
the expression on the right-hand side of the equality is equal to
$\langle l_2^*,[l,{\varepsilon}(l_1^*)]\rangle $.

\smallskip
{\bf LEMMA 4.} Let $l^*$ be from $L(A^*)$ and $l\in L(A)$. Then
there exists an element ${\rm i}\in I(A^*)$ satisfying the
equalities $[l,l^*]=[l,{\rm i}]$ and $\langle l^*$, $l\rangle
=\langle {\rm i}$, $l\rangle$.

{\bf Proof.} It suffices to prove the equalities for the case when
$l^*$ is in ${\rm WIntDer}(A^*)$. Let $l=\sum [f_i,a_i]$ and $g$
be an arbitrary element from $A^*$. Take the subcoalgebra $B$
in $\langle A,\bigtriangleup\rangle$ generated by elements $a_i$.
By Proposition 2, then, $Bl^*=B{\rm i}$ and $\langle
gl^*,B\rangle = \langle g{\rm i},B\rangle$ for some ${\rm i}$ in
${\rm IntDer}(A^*) $. Since $gl\in B$, we have $gll^*= gl{\rm
i}$. On the other hand, $gl^*l=\sum\limits_i (f_i,gl^*,a_i)$.
Let ${\Delta}(a_i)=\sum b^i_{(1)}\otimes b^i_{(2)}$, where
$b^i_{(1)}, \;b^i_{(2)} \in B$. Then
$$ \sum\limits_{i} (f_i,gl^*,a_i)=\sum\limits_i\sum (b^i_{(1)}
\langle f_i(gl^*), b^i_{(2)}\rangle - (f_i\cdot b^i_{(1)}) \langle
gl^*, b^i_{(2)}\rangle)= $$ $$ \sum\limits_i\sum (b^i_{(1)} \langle
gl^*,(f_i\cdot b^i_{(2)}) \rangle - (f_i\cdot b^i_{(1)})\langle
gl^*,b^i_{(2)}\rangle)= $$ $$ \sum\limits_i\sum (b^i_{(1)}\langle
g{\rm i}, (f_i\cdot b^i_{(2)})\rangle - (f_i\cdot b^i_{(1)})\langle
g{\rm i}, b^i_{(2)}\rangle)= \sum\limits_i (f_i, g{\rm i},a_i)=
gl{\rm i}. $$

\noindent Consequently, $g[l,l^*]=g[l,{\rm i}]$, thus proving the
first equality. The second equality is evident.

\medskip
\begin{center}
{\bf 2. THE KKT CONSTRUCTION FOR JORDAN ALGEBRAS}
\end{center}
\smallskip

Let $\rho$ be a standard embedding of the space $L(A^*)\otimes
L(A^*)$ into $(L(A)\otimes L(A))^*$. Consider the
finite-dimensional subcoalgebra $B$ of $\langle A,{\Delta}\rangle$
and the linear space $L=B\oplus B'\oplus [A^*,B]\oplus \bar{B}$.
Now if ${\rm d}$ is an element in ${\rm WIntDer}(A^*)$, then $A^*
({\rm d}-{\rm i})\subseteq B^\bot$ for some inner derivation ${\rm
i}$. Therefore, it is not hard to see that $L(A^*)=I(A^*)+L^\bot$.

\smallskip
{\bf LEMMA 5.} There exists a map ${\Delta}_B$ from $L$ into
$L(A)\otimes L(A)$ such that for any $f$, $g\in L(A^*)$ and
$l\in L$,
$$ \langle [f,g],l\rangle =\langle \rho (f\otimes g), {\Delta}_B(l) \rangle. $$

{\bf Proof.} Since $B$ is a finite-dimensional space, it follows
that $B^\bot$, the orthogonal complement of $B$ in $A^*$, has
finite codimension. Therefore, $A^*=V\oplus B^\bot$, where $V$
is a finite-dimensional subspace. The space $L$ has finite
dimension by Proposition 1. Obviously, $L(A)=L\oplus L_1$. It
follows from Lemma 4 that $L$ is an $L(A^*)$-subbimodule of
$L(A)$. By Lemma 3, $L^\bot$, the orthogonal complement of $L$,
is an ideal in $L(A^*)$. The ideal generated in $L(A^*)$ by the
set $B^\bot$ is obviously contained in $L^\bot$. Therefore,
$L(A^*)=V\oplus V'\oplus [V,V]\oplus\bar{V}\oplus L^\bot$. By
Proposition 5, we can suppose the space $V\oplus V'\oplus
[V,V]\oplus \bar{V}$ to be contained in $L^*+[A^*,B]^\bot$. And
it is now possible to assume that $L(A^*)= L^*\oplus L^\bot$.

Let $\rho (L(A^*)\otimes L(A^*))_0$ be the projection of the space
$\rho (L(A^*)\otimes L(A^*))$ onto $(L\otimes L)^*$ and let $f_L$
be the projection of an element $f$ from $L(A^*)$ onto the
subspace $L^*$. It is then easy to see that there exists an
injective map $\rho_B\!: L^*\otimes L^*\to (L\otimes L)^*$, in
which case $\rho_B(f_L\otimes g_L)=\rho (f\otimes g)_0$ for any $f$
and $g$ in $L(A^*)$. Consequently, the dual map $\rho_B^{*}\!:
L\otimes L\to (L^*\otimes L^*)^*$ is invertible. Consider the
linear map $M\!: L^*\otimes L^*\to L^*$, defined by $M(f_L\otimes
g_L)= [f,g]_L$. Since $L^\bot$ is an ideal in $L(A^*)$, and
$L^*\cap L^\bot=0$, the $M$ is well defined. Suppose $M^*$ is
the dual of the map $M$. Then $M^*$ is the map from $L$ into
$(L^*\otimes L^*)^*$, in which case $\langle f_L\otimes g_L $,
$M^*(n)\rangle =\langle [f,g]_L$, $n\rangle$ for any elements $f$,
$g\in L(A^*)$ and $n\in L$.

Put ${\Delta}_B=(\rho_B^*)^{-1}M^*$. Clearly, ${\Delta}_B$ is a
map from $L$ into $L\otimes L$, and $\langle [f,g]_L$,
$l\rangle =\langle \rho_B (f_L\otimes g_L)$,
${\Delta}_B(l)\rangle$ for any elements $f$, $g$ in $L(A^*)$
and $l\in L$. Consider elements $f$ and $g$ from $L(A^*)$.
Then $f=f_L+l_1$ and $g=g_L+l_2$, where $f_L$, $g_L\in L^*$,
$l_1$, $l_2\in L^\bot$. Since $\rho (L^*\otimes L^\bot
+L^\bot\otimes L^*+L^\bot \otimes L^\bot)\subseteq (L\otimes
L)^\bot$ and $L^\bot$ is an ideal in $L(A^*)$, for any elements
$f$, $g$ in $L(A^*)$ and $l$ in $L$ we obtain $\langle \rho
(f\otimes g)$, ${\Delta}_B (l)\rangle = \langle \rho (f\otimes
g)_0$, ${\Delta}_B(l)\rangle = \langle \rho_B (f_L\otimes
g_L)$, ${\Delta}_B(l)\rangle =\langle [f,g]_L$, $l\rangle=\langle
[f,g],l\rangle$.

We proceed to the main theorem.

\smallskip
{\bf THEOREM 1.} On the linear space $L(A)$, comultiplication
${\Delta}_L$ can be defined in such a way that the coalgebra
$\langle L(A), {\Delta}_L\rangle$ is a Lie coalgebra with the dual
algebra $L(A^*)$.

{\bf Proof.} Let $l\in L(A)$ and $l=a+b'+\sum[f_i,a_i]+ \bar{c}$,
where $a$, $b$, $c$, $a_i\in A$ and $f_i\in A^*$. Take the
subcoalgebra $B$ of $\langle A,{\Delta}\rangle$, generated by
elements $a$, $b$, $c$ and by all $a_i$. Then $B$ is a
finite-dimensional subspace. Let $L=B\oplus B'\oplus [A^*,B]\oplus
\bar{B}$. Obviously, $l\in L$. Define comultiplication
${\Delta}_L$ by setting ${\Delta}_L(l)={\Delta}_B(l)$, where
${\Delta}_B$ is the map specified in Lemma 5. The map ${\Delta}_L$
is well defined. Indeed, let $l_1\in L(A)$ and $l=l_1$. It is
then easy to see that $l_1= a+b'+\sum [g_i,e_i]+\bar{c}$, where
$g_i\in L(A^*)$ and $e_i\in L(A)$. Take the subcoalgebra $D$ of
$\langle A,{\Delta}\rangle $, generated by elements $a$, $b$,
$c$ and by all $e_i$. Let $L_1=D\oplus D'\oplus [A^*,D]\oplus
\bar{D}$. By Lemma 5, we obtain $\langle \rho (f\otimes g)$,
${\Delta}_B(l)\rangle =\langle [f,g]$, $l\rangle = \langle \rho
(f\otimes g)$, ${\Delta}_D(l)\rangle$ for any $f$, $g\in L(A^*)$.
Since the space $\rho(L(A^*)\otimes L(A^*))$ is dense in
$(L(A)\otimes L(A))^*$, the latter equality yields ${\Delta}_B(l)=
{\Delta}_D(l)$. This means that ${\Delta}_L$ is well defined, and by
Lemma 5, $\langle [f,g]$, $l\rangle =\langle \rho (f\otimes g),
{\Delta}_L(l)\rangle$ for any $f$, $g\in L(A^*)$ and $l\in L(A)$.
Now it is not hard to see that ${\Delta}_L$ is the desired
comultiplication.

\medskip
\begin{center}
{\bf 3. SUBCOALGEBRAS AND COIDEALS OF $\langle
L(A),{\Delta}_L\rangle$} \end{center} \smallskip

In this section we reason to establish the relationship between
coideals (subcoalgebras) of the Jordan algebra $\langle
A,{\Delta}\rangle$ and the coalgebra $\langle
L(A),{\Delta}_L\rangle$.

Let $B$ be a subalgebra of $A^*$ such that $1\in B$. Then, as is
known (see [4]), the spaces $L_1(B)=B\oplus B'\oplus [B,B]\oplus
\bar{B}$ and $L_2(B)= B\oplus B'\oplus {\rm d}(B)\oplus \bar{B}$,
where ${\rm d}(B)=\{{\rm d}\in {\rm WIntDer}(A^*)|$ $B{\rm
d}\subseteq B\}$, are subalgebras in $L(A^*)$. Either
subalgebra contains a three-dimensional simple Lie algebra $U$
generated in $L(A^*)$ by the elements 1, $(1)'$, and $\bar{1}$.
Conversely, if $L$ is a subalgebra in $L(A^*)$ that is closed
under the map ${\varepsilon}$ and contains a subalgebra
$U$, then $B=A^*\cap L$ is a subalgebra in $A^*$, $1\in B$,
and $L_1(B)\subseteq L\subseteq L_2 (B)$. Obviously, $L_1(B)$ is
an ideal in $L_2(B)$.

First let $V$ be a coideal of the coalgebra $\langle
A,{\Delta}\rangle$. Then the orthogonal complement $V^\bot$ of
the space $V$ in $A^*$ is a subalgebra, and $1\in V^\bot$. It
is not hard to note that $V^\bot\cdot V\subseteq V$. Consequently,
$V^\bot [V^\bot, V]\subseteq V$. Consider the spaces $L_1^*(V)=
V\oplus V'\oplus [V^\bot, V]\oplus \bar{V}$ and $L_2^*(V)= V\oplus
V' \oplus {\rm d}^* (V)\oplus \bar{V}$, where ${\rm d}^*(V)=\{
v\in [A^*,A]|$ $V^\bot v\subseteq V\}$. Obviously,
$L_1^*(V)\subseteq L_2^*(V)$.

We show that $L_1^*(V)^\bot=L_2(V^\bot )$ and $L_2^*(V)^\bot =
L_1(V^\bot)$. Indeed, by Lemma 3, $\langle {\rm d}(V^\bot)$,
$[V^\bot,V]\rangle \subseteq \langle [V^\bot,{\rm d}(V^\bot)]$,
$V'\rangle \subseteq \langle V^\bot {\rm d}(V^\bot)$, $V\rangle
\subseteq \langle V^\bot$, $V\rangle =0$. Therefore, $[V^\bot,V]
\subseteq {\rm d}(V^\bot)^\bot$. Again, if ${\rm d}\in{\rm
WIntDer}(A^*)$ and ${\rm d}\in [V^\bot,V]^\bot$, by Lemma 3 we
have $\langle V^\bot{\rm d}$, $V\rangle \subseteq \langle [V^\bot,
{\rm d}]$, $V\rangle \subseteq \langle {\rm d}$, $[V^\bot,
V]\rangle =0$. Therefore, $V^\bot {\rm d}\subseteq V^\bot$, that
is, ${\rm d}\in {\rm d}(V^\bot )$. This implies that
$[V^\bot,V]^\bot\subseteq {\rm d}(V^\bot)$, where
$[V^\bot,V]^\bot$ is an orthogonal complement of the space
$[V^\bot,V]$ in ${\rm WIntDer}(A^*)$. Consequently, $[V^\bot,V]=
{\rm d}(V^\bot)^\bot$, and it is now easy to see that $L_1^*
(V)^\bot= L_2 (V^\bot)$. The second equality is proved similarly.

Since $L_1(V^\bot)$ and $L_2(V^\bot )$ are subalgebras of the Lie
algebra $L(A^*)$, $L_1^*(V)$ and $L_2^*(V)$ are coideals in
$L(A)$. The coideals are obviously closed under the map $\pi$.
Therefore, the coalgebras $L_1(V^\bot)$ and $L_2(V^\bot)$ are
closed under ${\varepsilon}$, and since $1\in V^\bot$, the
subalgebra $U$ is contained in $L_1(V^\bot)$ and $L_2(V^\bot)$.
Consequently, the coideals $L_1^*(V)$ and $L_2^*(V)$ are contained
in $U^\bot$.

Next let $L$ be a coideal in $\langle L(A),{\Delta}_L\rangle$ which
is closed under the map ${\pi}$ and is such that $L\subseteq
U^\bot$. Then $L^\bot$ is a subalgebra in $L(A)$ that is closed
with respect to ${\varepsilon}$, and $U\subseteq L^\bot$.
Consequently, for the subalgebra $B=A^*\cap L^\bot$ we have
$L_1(B)\subseteq L^\bot\subseteq L_2(B)$. Since $B=V^\bot$ for
some coideal $V$ of $\langle A$, ${\Delta}\rangle$, we obtain
$L_1^*(V)\subseteq L\subseteq L_2^*(V)$ by the above argument.
Obviously, $A\cap L=V$. Thus, the following is valid:

\smallskip
{\bf THEOREM 2.} Let $V$ be a coideal of the coalgebra $\langle
A,{\Delta}\rangle$. Then the spaces $L_1^*(V)$ and $L_2^*(V)$
are coideals of the coalgebra $\langle L(A)$, ${\Delta}_L\rangle$,
in which case also $L_1^*(V)^\bot=L_2(V^\bot)$ and
$L_2^*(V)^\bot=L_1(V^\bot)$. Conversely, if $L$ is a coideal in
$\langle L(A)$, ${\Delta}_L\rangle$ that is closed under ${\pi}$
and is such that $L\subseteq U^\bot $, then $V=A\cap L$ is a
coideal in $\langle A$, ${\Delta}\rangle$, and
$L_1^*(V)\subseteq L\subseteq L_2^*(V)$.

Assume that $B$ is an ideal of the algebra $A^*$. Take the set $p(B)=
\{{\rm d}\in {\rm WIntDer}(A^*)| (A^*) {\rm d}\subseteq B\}$. We
show that $p(B)$ is an ideal in ${\rm WIntDer}(A^*)$. In fact,
let ${\rm d}$ and ${\rm d}_1$ be arbitrary elements from $p(B)$
and ${\rm WIntDer}(A)$, respectively. By Proposition 2, we have
$(B^\bot){\rm d}_1\subseteq B^\bot$ and $\langle B{\rm d}_1$,
$B^\bot\rangle \subseteq \langle B$, $(B^\bot){\rm d}_1\rangle =0$.
Therefore, $\langle (A^*)[{\rm d},{\rm d}_1]$, $B^\bot\rangle
\subseteq \langle B{\rm d}_1$, $B^\bot\rangle =0$. This implies
that $(A^*)[{\rm d},{\rm d}_1]\subseteq B$, that is, $p(B)$ is an
ideal in ${\rm WIntDer}(A^*)$. It is then not difficult to see
that the spaces $I_1(B)=B\oplus B'\oplus [A^*,B]\oplus \bar{B}$ and
$I_2(B)= B\oplus B'\oplus p(B)\oplus \bar{B}$ are ideals in
$L(A^*)$. Conversely, if $I$ is an ideal in $L(A^*)$, and $B=
A^*\cap I$, appealing to the argument in [4, p. 331], it is not
hard to prove that $B$ is an ideal of $A^*$, and $I_1(B)\subseteq
I\subseteq I_2(B)$.

We turn to a subcoalgebra $B$ of the coalgebra $\langle
A,{\Delta}\rangle$. Clearly, the orthogonal complement $B^\bot$
of the space $B$ in $A^*$ is an ideal in $A^*$. By Proposition
2, $B{\rm d}\subseteq B$ for any ${\rm d}$ from ${\rm
WIntDer}(A^*)$. Consider the spaces
$$ I_1^*(B)=B\oplus B'\oplus [A^*,B]\oplus \bar{B}\mbox{\, and \,}
I_2^*(B)=B\oplus B'\oplus p^*(B)\oplus \bar{B}, $$

\noindent where $p^*(B)=\{v\in [A^*,A]|$ $A^*v\subseteq B\}$. It
is not hard to see that $p^*(B)$ is a ${\rm WIntDer}(A)$-subbimodule
of the Lie ${\rm WIntDer}(A^*)$-bimodule $[A^*,A]$.
It is also evident that $I_1^*(B)\subseteq I_2^*(B)$.

We show that $I_1^*(B)^\bot=I_2(B^\bot)$ and $I_2^*(B)^\bot=
I_1(B^\bot)$. Indeed, by Lemma 3, $\langle p(B^\bot)$,
$[A^*,B]\rangle \subseteq \langle [A^*,p(B^\bot)]$,
$B'\rangle\subseteq\langle A^*p(B^\bot)$,
$B\rangle\subseteq\langle B^\bot$, $B\rangle =0$. Therefore,
$[A^*,B]\subseteq p(B^\bot)^\bot$. Again, if ${\rm d}\in {\rm
WIntDer}(A^*)$ and ${\rm d}\in [A^*,B]^\bot$, by Lemma 3 we have
$\langle A^*{\rm d}$, $B\rangle\subseteq \langle [A^*,{\rm d}]$,
$B\rangle\subseteq \langle {\rm d}$, $[A^*,B]\rangle =0$. Hence,
$A^*{\rm d}\subseteq B^\bot$, that is, ${\rm d}\in p(B^\bot )$.
This implies that $[A^*,B]^\bot\subseteq p(B^\bot)$, where
$[A^*,B]^\bot$ is an orthogonal complement of the space $[A^*,B]$
in ${\rm WIntDer}(A^*)$. Consequently, $[A^*,B]^\bot=p(B^\bot)$,
and then it is easy to see that $I_1^*(B)^\bot=I_2(B^\bot )$.
The second equality is proved similarly. Since $I_1(B^\bot)$ and
$I_2(B^\bot)$ are ideals of the Lie algebra $L(A^*)$, it follows
that $I_1^*(B)$ and $I_2^*(B)$ are subcoalgebras in $L(A)$.

Now let $I$ be a subcoalgebra of $\langle L(A), {\Delta}_L \rangle$.
Then $I^\bot$ is an ideal in $L(A^*)$. Consequently, for the
ideal $B=A^*\cap I^\bot$ in $A^*$ we have $I_1(B)\subseteq
I^\bot\subseteq I_2(B)$. Since $B=V^\bot$ for some subcoalgebra
$V$ of $\langle A,{\Delta}\rangle $, the above argument implies
$I_1^*(V)\subseteq I\subseteq I_2^*(V)$. Obviously, $A\cap I=V$.
The following is thus valid:

\smallskip
{\bf THEOREM 3.} Let $B$ be a subcoalgebra of the coalgebra
$\langle A,{\Delta}\rangle $. Then the spaces $I_1^*(B)$ and
$I_2^*(B)$ are subcoalgebras of $\langle L(A)$,
${\Delta}_L\rangle$, in which case also $I_1^*(B)^\bot=
I_2(B^\bot)$ and $I_2^*(B)^\bot=I_1(B^\bot )$. Conversely, if $I$
is a subcoalgebra in $\langle L(A)$, ${\Delta}_L\rangle$, then
$B=A\cap I$ is a subcoalgebra of $A$, and $I_1^*(B)\subseteq
I\subseteq I_2^*(B)$.

\medskip
\begin{center}
{\bf REFERENCES}
\end{center}
\begin{itemize}

\item V. G. Drinfeld, ``Hamiltonian structures on Lie groups, Lie
bialgebras, and geometric meaning of the classical Yang--Baxter
equations," {\it Dokl. Akad. Nauk SSSR}, {\bf 268}, No. 2, 285-287
(1983).

\item J. Anquela, T. Cortes, and F. Montaner, ``Nonassociative
coalgebras,'' {\it Comm. Alg.}, {\bf 22}, No. 12, 4693-4716 (1994).

\item M. E. Sweedler, {\it Hopf Algebras}, Benjamin, New York
(1969).

\item N. Jacobson, ``Structure and representations of Jordan
algebras,'' in {\it Am. Math. Soc. Coll. Publ.,} {\bf 39},
Providence, Rhode Island (1968).
\end{itemize}

\end{document}